\documentstyle[amsfonts,12pt]{article}
\input amssym.def

\newcommand{\End}{{\rm E n d }}

\newcommand{\Mat}{{\rm M a t }}

\newcommand{\ES}{{\rm E}{\cal S}}

\newcommand{\brarrow}{\succ\rightarrow}

\newcommand{\bbrarrow}{\succ\succ\rightarrow}
\newcommand{\bblarrow}{\leftarrow\prec\prec}

\newcommand{\SM}{{\cal S}M}

\newcommand{\Tp}{T_{poly}^{\bul}}

\newcommand{\cTp}{\cT^{\bul}_{poly}}

\newcommand{\cAb}{{\cal A}^{\bullet}}
\newcommand{\Omb}{\Om^{\bul}}

\newcommand{\OmT}{\Om^{\bul}(\cT_{poly})}
\newcommand{\OmD}{\Om^{\bul}(\Cbu(\SM))}
\newcommand{\OmC}{\Om^{\bul}(\Cbd(\SM))}
\newcommand{\OmE}{\Om^{\bul}(\cE)}
\newcommand{\OM}{\cO_M}

\newcommand{\Cbu}{C^{\bullet}}
\newcommand{\Cbd}{C_{\bullet}}

\newcommand{\Linf}{L_{\infty}}

\newcommand{\tga}{\widetilde{\gamma}}

\newcommand{\tD}{\widetilde{D}}

\newcommand{\al}{{\alpha}}
\newcommand{\la}{{\lambda}}

\newcommand{\bul}{{\bullet}}

\newcommand{\mb}{{\mathfrak{b}}}

\newcommand{\Om}{{\Omega}}

\newcommand{\ga}{{\gamma}}

\newcommand{\G}{{\Gamma}}

\newcommand{\pa}{{\partial}}
\newcommand{\tr}{{\rm t r}}

\newcommand{\cK}{{\cal K}}
\newcommand{\cC}{{\cal C}}
\newcommand{\cM}{{\cal M}}

\newcommand{\cL}{{\cal L}}
\newcommand{\cD}{{\cal D}}
\newcommand{\cA}{{\cal A}}
\newcommand{\cG}{{\cal G}}

\newcommand{\cT}{{\cal T}}

\newcommand{\cS}{{\cal S}}
\newcommand{\cE}{{\cal E}}

\newcommand{\cO}{{\cal O}}

\newcommand{\bbR}{{\Bbb R}}

\newcommand{\bbA}{{\Bbb A}}

\newcommand{\n}{{\nabla}}

\newcommand{\de}{{\delta}}
\newcommand{\D}{{\Delta}}

\newcommand{\trd}{{\rm t r d}}
\newcommand{\ind}{{\rm i n d}}
\newcommand{\cotr}{{\rm c o t r}}

\date{}

\newtheorem{lem}{Lemma}
\newtheorem{teo}{Theorem}

\newtheorem{pred}{Proposition}

\title{Formality theorem for Hochschild
(co)chains of the algebra of endomorphisms
of a vector bundle}

\author{Vasiliy Dolgushev}

\begin{document}

\maketitle

\begin{abstract}
We prove the formality theorem for the differential
graded Lie algebra module of Hochschild chains 
for the algebra of endomorphisms of a smooth 
vector bundle. We discuss a possible application 
of this result to a version of the algebraic index 
theorem for Poisson manifolds.
\end{abstract}
~\\[0.3cm]

\section{Introduction}
The purpose of this note is to prove
that the differential
graded Lie algebra (DGLA) module of Hochschild chains 
for the algebra of endomorphisms of a smooth 
vector bundle over a manifold $M$ can be 
connected by a chain of quasi-isomorphisms to its 
cohomology. 

For a trivial vector bundle one
can easily prove the desired statement using
the formality theorem \cite{FTHC},
\cite{thesis} for Hochschild
chains of the algebra $\cO_M$ of functions
on $M$ and the (co)trace map \cite{Loday} between
the Hochschild complexes of the algebra $\OM$
and the algebra $\Mat(\OM)$ of finite size
matrices with entries in $\OM$\,.

However, for a non-trivial vector bundle $E$
both the trace and the cotrace maps
between the Hochschild complexes of
the algebra of functions and the algebra 
of endomorphisms
are defined only locally. Thus, the
question of formality of the 
DGLA module of Hochschild chains for the 
algebra of endomorphisms of a vector bundle
requires some work.

In~ this~ note~ we~ propose~ Fedosov's~ resolution~
\cite{thesis}~
of~ the~
sheaf~ of~ DGLA~ modules \\ $(\Cbu(\End(E)),\Cbd(\End(E)))$
of Hochschild (co)chains for the sheaf $\End(E)$ of endomorphisms
of a vector bundle $E$ and
construct a quasi-isomorphism between this resolution and
Fedosov's resolution of the sheaf of DGLA modules
$(\Cbu(\OM),\Cbd(\OM))$ for $\OM$\,. Combining this construction
with the result in \cite{thesis} (see diagram
(5.15) on page 86) we obtain a proof of the
formality theorem for the sheaf of
DGLA modules $(\Cbu(\End(E)),\Cbd(\End(E)))$\,.

Besides the obvious applications  
of this theorem to the questions of 
deformation quantization this result may lead to
an interesting version of the algebraic
index theorem for Poisson manifolds.
We will devote a separate paper to this 
version of the algebraic index theorem. 
Here we will only outline a rough idea. 

There are many versions of the algebraic
index theorem for symplectic manifolds
\cite{BNT}, \cite{BNT1}, \cite{Ch-I}, \cite{Fedosov}
\cite{NT}. One of these versions 
\cite{Ch-I} describes a natural map 
(see eq. (31) and theorem 4 in \cite{Ch-I})
$$
c l : K_0(\bbA) \to H^{top}_{DR}(M)
$$ 
from the $K$-theory of the deformation quantization 
algebra $\bbA$ to the top degree De Rham cohomology of $M$\,.
This map is obtained by composing 
the trace density map \cite{FFS} 
\begin{equation}
\label{trd}
\trd : \bbA/ [\bbA, \bbA]  \to  H^{top}_{DR}(M)
\end{equation}
from the zeroth Hochschild homology 
$HH_0(\bbA) = \bbA/ [\bbA, \bbA]$ of $\bbA$ to 
the top degree De Rham cohomology of $M$ with
the lowest component of the Chern character
(see example 8.3.6 in \cite{Loday})
\begin{equation}
\label{Ch}
c h_{0,0} : K_0(\bbA) \to \bbA/ [\bbA, \bbA] 
\end{equation}
from the $K$-theory of $\bbA$ to the 
zeroth Hochschild homology of $\bbA$\,. 

In the case of a general Poisson (not symplectic)
manifold $M$ one cannot construct 
the trace density map (\ref{trd}). 
Instead we have the map
\begin{equation}
\label{trd1}
\mu : \bbA/ [\bbA, \bbA]  \to  HP_0(M)
\end{equation}
from zeroth Hochschild homology of $\bbA$ 
to the zeroth Poisson homology \cite{JLB}, \cite{Koszul} 
of $M$. This map is constructed with the 
help\footnote{See corollary 2 on page 92 in \cite{thesis}.}
of the formality theorem for Hochschild chains 
\cite{thesis}, \cite{Sh} of $\OM$\,. 

Composing (\ref{trd1}) with (\ref{Ch}) we 
get the map
\begin{equation}
\label{ind}
\ind : K_0 (\bbA) \to HP_0(M) 
\end{equation}
from the $K$-theory of $\bbA$ to the 
zeroth Poisson homology of $M$\,.

Following the arguments of the proof of 
theorem 6.1.3 in \cite{Fedosov} one can 
show that the image $\ind(Q)$ of a $K$-theory
element $Q\in K_0(\bbA)$ depends only 
on the principal part of $Q$. The desired 
algebraic index theorem should give an 
explicit form of this dependence and  
we expect that the formality theorem 
for Hochschild chains of $\End(E)$ 
will help us to solve this problem.

The paper is organized as follows. 
In the second section we fix notation
and recall some results and constructions 
we are using here. In section 3 we formulate 
the main result of this paper (see theorem \ref{Main})
and give a proof. In the concluding section 
we propose an interesting generalization 
of lemma \ref{ona} which is used in the proof
of theorem \ref{Main}.  
We suspect that this generalization may shed 
some light on an unknown index formula 
for the map (\ref{ind})\,.

~\\
{\bf Acknowledgment.}
I would like to thank Volodya Rubtsov for 
stimulating discussions of this topic.  
I started this project when I was 
a Liftoff Fellow of Clay Mathematics Institute. 
I thank this Institute for the 
support. I am also partially supported by
the Grant for Support of Scientific
Schools NSh-8065.2006.2.

\section{Preliminaries}
In this section we fix
notation and recall some results we are going to
use in this note.

For an associative algebra $B$
we denote by $\Mat(B)$
the algebra of finite size matrices over
$B$\,. $\Cbd(B)$ is the normalized
Hochschild chain complex
of $B$ with coefficients in $B$
\begin{equation}
\label{chains}
\Cbd(B) = \Cbd(B,B)
\end{equation}
and $\Cbu(B)$ is the normalized
Hochschild cochain complex of $B$
with coefficients in $B$ and with shifted
grading
\begin{equation}
\label{rule}
\Cbu(B) = C^{\bul + 1}(B,B)\,.
\end{equation}
The Hochschild coboundary operator
is denoted by $\pa$ and the
Hochschild boundary operator is
denoted by $\mb$\,.

It is well known that the 
Hochschild cochain complex (\ref{rule}) 
carries the structure of a differential graded 
Lie algebra (DGLA). The corresponding 
Lie bracket (see eq. (3.2) on page 45 in \cite{thesis}) 
was originally introduced by M. Gerstenhaber 
in \cite{Ger}. We will denote this bracket by 
$[,]_G$\,. 

The Hochschild chain complex (\ref{chains}) 
carries the structure of a differential graded 
Lie algebra module over 
the DGLA $\Cbu(B)$\,. We will denote the 
action (see eq. (3.5) on page 46 in \cite{thesis}) 
of cochains on chains by $R$\,.

The trace map $\tr$ \cite{Loday} is the map from the
Hochschild chain complex $\Cbd(\Mat(B))$ of
the algebra $\Mat(B)$ to the Hochschild chain
complex $\Cbd(B)$ of the algebra $B$\,. This map
is defined by the formula
\begin{equation}
\label{trace-map}
\tr (M_0 \otimes M_1 \otimes \dots \otimes M_k) =
\sum_{i_0, \dots, i_k}
(M_0)_{i_0 i_1} \otimes (M_1)_{i_1 i_2} \otimes
\dots \otimes (M_k)_{i_k i_0}\,,
\end{equation}
where $M_0, \dots, M_k$ are matrices in $\Mat(B)$
and $(M_a)_{ij}$ are the corresponding entries.

Dually, the cotrace map \cite{Loday}
$$
\cotr : \Cbu(B) \to \Cbu(\Mat(B))
$$
is defined by the formula
\begin{equation}
\label{cotrace-map}
(\cotr(P)(M_0, M_1, \dots , M_k))_{ij} =
\sum_{i_1, \dots, i_k}
P((M_0)_{i i_1}, (M_1)_{i_1 i_2},
\dots , (M_k)_{i_k  j})\,,
\end{equation}
where $P\in C^{k}(B)$ and
$M_0, \dots, M_k$ are, as above, matrices
in $\Mat(B)$\,.

``DGLA'' always means a differential graded Lie algebra.
The arrow $\brarrow$ denotes an $\Linf$-morphism
of sheaves of DGLAs, 
the arrow $\bbrarrow$ denotes a
morphism of sheaves of $\Linf$-modules,  
and the notation
$$
\begin{array}{c}
\cL\\[0.3cm]
\downarrow_{\,mod}\\[0.3cm]
\cM
\end{array}
$$
means that $\cM$ is a sheaf of DGLA modules
over the sheaf of DGLAs $\cL$\,. The symbol
$\circ$ always stands for the composition of morphisms.

Throughout this note
$M$ is a smooth real ma\-ni\-fold
of dimension $d$\,. $E$ is a smooth real 
vector bundle over $M$ and $\End(E)$ denotes 
the sheaf of endomorphisms of $E$\,.
We denote by $\G(M, \cG)$  the vector space 
of sections of the sheaf
$\cG$ and by  $\Omb(\cG)$ the sheaf
of exterior forms with values in $\cG$\,.
We omit the symbol
$\wedge$ referring to a local basis of exterior forms,
as if one
thought of $dx^i$'s as anti-commuting variables.

$T^{\bul}_{poly}$ is the sheaf of polyvector
fields with shifted grading
$$
T^{\bul}_{poly} = \wedge^{\bul+1}_{\OM} TM\,,
\qquad  
T^{-1}_{poly} = \OM
$$ 
and $\cAb$ is the sheaf of exterior forms.   

$\Tp$ is a sheaf of graded Lie algebras 
with respect the so-called Schouten-Nijenhuis 
bracket $[,]_{SN}$
(see eq. (3.20) on page 50 in \cite{thesis}) and 
$\cAb$ is the sheaf of graded Lie algebra modules 
over $\Tp$ with respect to Lie derivative 
(see eq. (3.21) on page 51
in \cite{thesis}). We will regard $\Tp$ 
(resp. $\cAb$) as the sheaf of DGLAs 
(resp. the sheaf of DGLA modules)
with vanishing differential. 

We denote by $x^i$ local coordinates on $M$ and
by $y^i$ fiber coordinates in the tangent bundle
$TM$. Having these coordinates $y^i$ we
can introduce another local basis of exterior
forms $\{ dy^i \}$. We will use both bases
$ \{ dx^i \}$ and $\{ dy^i \}$. In particular,
the notation $\Omb(\cG)$
is  reserved for the sheaf of
$d y$-exterior forms with values in
the sheaf $\cG$ while $\cAb$ denotes the 
sheaf of $d x$-exterior forms.

$\SM$ is the formally completed symmetric algebra
of the cotangent bundle $T^*(M)$\,. Sections of the
sheaf $\SM$ can be viewed as formal power series in
tangent coordinates $y^i$\,. We regard $\SM$ as 
the sheaf of algebras over $\OM$\,. In particular, 
$\Cbu(\SM)$ is the sheaf of normalized 
Hochschild cochains of $\SM$ over $\OM$. Namely, 
the sections of $C^k(\SM)$ over an open subset 
$U\subset M$ are $\OM$-linear polydifferential 
operators with respect to tangent 
coordinates $y^i$
$$
P : \G(U,\SM)^{\otimes \, (k+1)} \to \G(U,\SM)
$$    
satisfying the normalization condition 
$$
P(\dots, f, \dots) = 0\,, 
\qquad \forall ~ f \in \OM(U)\,. 
$$
Similarly, $\Cbd(\SM)$ is the sheaf of normalized 
Hochschild chains\footnote{In \cite{thesis} the sheaf 
$\Cbu(\SM)$ is denoted by $\cD^{\bul}_{poly}$ and 
the sheaf $\Cbd(\SM)$ is denoted by $\cC^{poly}_{\bul}$.} 
of $\SM$ over $\OM$\,.
As in \cite{thesis} the tensor product in 
$$
C_k(\SM) = 
\underbrace{\SM \hat{\otimes}_{\OM} (\SM / \OM) 
\hat{\otimes}_{\OM}
\dots 
\hat{\otimes}_{\OM} (\SM / \OM)  }_{k+1}
$$
is completed in the adic topology in fiber 
coordinates $y^i$ on the tangent bundle $TM$\,.

The cohomology of the complex of sheaves 
$\Cbu(\SM)$ is the sheaf $\cTp$ of fiberwise 
polyvector fields (see page 60 in \cite{thesis}). 
The cohomology of the complex of sheaves
$\Cbd(\SM)$ is the sheaf $\cE$ of fiberwise differential 
forms (see page 62 in \cite{thesis}). These are 
$d x$-forms with values in $\SM$\,. 

In \cite{thesis} (see theorem 4 on page 68) it 
is shown that the sheaf of algebras $\Omb(\SM)$
can be equipped with a differential of the 
following form
\begin{equation}
\label{DDD}
D = \n - \de + A\,,
\end{equation}
where
\begin{equation}
\label{nabla}
\n = dy^i \frac{\pa}{\pa x^i} -
dy^i \G^k_{ij}(x) y^j \frac{\pa}{\pa y^k}\,,
\end{equation}
is a torsion free connection with
Christoffel symbols $\G^k_{ij}(x)$,
\begin{equation}
\label{delta}
\de = dy^i \frac{\pa}{\pa y^i}\,,
\end{equation}
and
$$
A=\sum_{p=2}^{\infty}dy^k A^j_{ki_1\dots i_p}(x)
y^{i_1} \dots
y^{i_p}\frac{\pa}{\pa y^j} \in \Om^1(M,\cT^0_{poly})\,.
$$
We refer to (\ref{DDD}) as the Fedosov differential.  

Notice that $\de$ in (\ref{delta}) is also a differential on 
$\Omb(\SM)$ and (\ref{DDD}) can be viewed as 
deformation of $\de$ via the connection $\n$\,. 

Let us recall from \cite{thesis}
the following operator 
on\footnote{The arrow over $\pa$ in (\ref{del-1})
means that we use the left derivative with respect to
the anti-com\-muting variable $d y^k$.} $\Omb(\SM)$
\begin{equation}
\de^{-1}(a) =
\cases{
\begin{array}{cc}
\displaystyle
y^k \frac {\vec{\partial}} {\partial (d y^k)}
\int\limits_0^1 a(x,t y,t d y)\frac{d t} t, & {\rm if}~
a \in \Om^{>0}(U, \SM)\,,\\
0, & {\rm otherwise}\,,
\end{array}
}
\label{del-1}
\end{equation}
which is used to prove the 
acyclicity of $\de$ and $D$ in positive 
dimension.

According to proposition $10$ on
page $64$ in \cite{thesis} the sheaves $\cTp$, $\Cbu(\SM)$,
$\cE^{\bul}$, and $\Cbd(\SM)$ are equipped with the
canonical action of the sheaf of Lie algebras
$\cT_{poly}^0$ and this action is compatible
with the corresponding (DG) algebraic structures.
Using this action in
chapter $4$ of \cite{thesis} we extend
the Fedosov differential (\ref{DDD}) to 
a differential on the sheaves 
$\Omb(\cTp)$, $\Omb(\cE^{\bul})$, 
$\Omb(\Cbd(\SM))$, and $\Omb(\Cbu(\SM))$
of DGLAs (resp. DGLA modules).

Using acyclicity of the Fedosov differential (\ref{DDD})
in positive dimension 
one constructs in \cite{thesis} embeddings
of the sheaves of 
DGLA modules\footnote{See eq. (5.1) on page 81 in \cite{thesis}.}
\begin{equation}
\begin{array}{ccc}
\Tp &\stackrel{\la_T}{\,\longrightarrow\,} &
(\OmT, D, [,]_{SN})\\[0.3cm]
\downarrow^{L}_{\,mod}  & ~  &
\downarrow^{L}_{\,mod}\\[0.3cm]
\cAb  &\stackrel{\la_{\cA}}{\,\longrightarrow\,} &
(\OmE, D),
\end{array}
\label{diag-T}
\end{equation}

\begin{equation}
\begin{array}{ccc}
(\OmD, D+\pa, [,]_{G}) &\stackrel{\,\la_D}{\,\longleftarrow\,} 
& \Cbu(\OM)\\[0.3cm]
\downarrow^{R}_{\,mod}  & ~  &     \downarrow^{R}_{\,mod} \\[0.3cm]
(\OmC, D+\mb) &\stackrel{\,\la_C}{\,\longleftarrow\,} &  \Cbd(\OM),
\end{array}
\label{diag-D}
\end{equation}
and shows that these are quasi-isomorphisms of the
corresponding complexes of sheaves.

Furthermore, using Kontsevich's and Shoikhet's
formality theorems for $\bbR^d$ \cite{K}, \cite{Sh}
in \cite{thesis} one constructs the following
diagram
\begin{equation}
\label{diag-K-Sh}
\begin{array}{ccc}
(\OmT, D, [,]_{SN}) & \stackrel{\cK}{\brarrow} &
(\OmD, D+\pa, [,]_{G}) \\[0.3cm]
\downarrow^{L}_{\,mod}  & ~  &
\downarrow^{R}_{\,mod} \\[0.3cm]
(\OmE, D) & \stackrel{\cS}{\bblarrow} &
(\OmC, D+\mb)
\end{array}
\end{equation}
where $\cK$ is an $\Linf$ quasi-isomorphism of
sheaves of DGLAs and $\cS$ is a quasi-isomorphism
of sheaves of $\Linf$-modules  over the
sheaf of DGLAs $(\OmT, D, [,]_{SN})$\,,
and the $\Linf$-module on $\OmC$ is obtained
by composing the quasi-isomorphism
$\cK$ with the DGLA modules structure $R$
(see eq. (3.5) on p. 46 in \cite{thesis} for the 
definition of $R$)\,.

Diagrams (\ref{diag-T}), (\ref{diag-D}) and (\ref{diag-K-Sh})
show that the sheaf $\Cbd(\OM)$ of DGLA modules of Hochschild
chains of $\OM$ is quasi-isomorphic to the
sheaf of graded Lie algebra modules $\cAb$ of its
cohomology.

~\\
{\bf Remark 1.} As in \cite{thesis} we use adapted versions 
of Hochschild (co)chains for the sheaves $\OM$ and $\End(E)$
of functions and of endomorphisms of a vector bundle $E$, 
respectively. Thus, $\Cbu(\OM)$ is a sheaf of polydifferential 
operators (see page 48 in \cite{thesis}) satisfying 
the corresponding normalization condition.
$\Cbu(\End(E))$ is the sheaf of (normalized) 
polydifferential operators acting on $\End(E)$ 
with coefficients in $\End(E)$\,.
Furthermore,  $\Cbd(\OM)$ is the sheaf of
(normalized) polyjets 
$$
C_k(\OM) = Hom_{\OM} (C^{k-1}(\OM), \OM)\,,
$$
and 
$$
C_k(\End(E)) = Hom_{\End(E)} (C^{k-1}(\End(E)), \End(E))\,.
$$
We have to warn the reader that the space of 
global sections of the sheaf $\Cbu(\OM)$ 
(resp. $\Cbu(\End(E))$) is not isomorphic to the 
space of Hochschild cochains of the algebra 
$\OM(M)$ of functions (resp. the algebra 
$\G(M, \End(E))$ of endomorphisms of $E$). 
Similar expectation is wrong for Hochschild 
chains. Instead we have the following inclusions: 
$$
\G(M, \Cbu(\OM)) \subset \Cbu(\OM(M))\,, 
\qquad 
\G(M, \Cbu(\End(E) )) \subset \Cbu(\G(M,\End(E)))\,, 
$$
$$
\Cbd(\OM(M)) \subset \G(M, \Cbd(\OM))\,, 
\qquad
\Cbd(\G(M,\End(E)) ) \subset \G(M, \Cbd(\End(E) ))\,. 
$$

~\\
{\bf Remark 2.} Unlike in \cite{thesis} we use 
only normalized Hochschild (co)chains. It is not hard 
to check that the results we need from \cite{thesis}, 
\cite{K}, and \cite{Sh} 
also hold when this normalization condition is imposed.

\section{The formality theorem}
Let $E$ be a smooth real vector bundle over
the smooth real manifold $M$ and let
$\End(E)$ denote the sheaf of endomorphisms
of $E$\,. We regard $\End(E)$ as
a sheaf of algebras over $\bbR$.

Here is the main result of this note: 
\begin{teo}
\label{Main}
The sheaf of DGLA modules $\Cbd(\End(E))$ over
the sheaf of DGLAs $\Cbu(\End(E))$ is formal.
\end{teo}
{\bf Proof.}
Let us introduce the following auxiliary sheaf of
algebras
\begin{equation}
\label{ES}
\ES = \End(E)\otimes_{\OM} \SM \,.
\end{equation}
Regarding $\ES$ as a sheaf of algebras over 
$\OM$ we also consider the following adapted versions 
of (normalized) Hochschild (co)chains. Thus, $C^k(\ES)$ is the 
sheaf whose sections over an open subset $U\subset M$
are $\OM$-polylinear maps 
$$
P: \G(U,\ES)^{\otimes\, {k+1}} \to \G(U,\ES)\,,
$$
which are differential in fiber coordinates $y^i$
and satisfy the normalization condition: 
$$
P (\dots, f, \dots ) = 0\,, \qquad 
\forall \quad f \in \OM(U)\,.
$$ 

Similarly, $C_k(\ES)$ is the sheaf of normalized 
Hochschild chains of $\ES$ over $\OM$ 
for which the tensor product is completed in 
the adic topology in fiber coordinates $y^i$ on 
the tangent bundle $TM$\,.

It is clear that the differentials $\pa$, $\mb$ 
as well as the operations $[,]_G$ and $R$ 
are well defined on the sheaves 
$\Cbu(\ES)$, $\Cbd(\ES)$. Thus, we regard $\Cbd(\ES)$ 
as the sheaf of DGLA modules over the sheaf of
DGLAs $\Cbu(\ES)$\,.

It is not hard to show that one can extend
the Fedosov differential (\ref{DDD}) on $\SM$ to
a differential on $\ES$ in the framework of
the following ansatz:
\begin{equation}
\label{DDD-E1}
\tD = D + [\ga^E,\, ]\,, \qquad
\ga^E= \G^E + \tga^E\,,
\end{equation}
where $\G^E$ is a connection form of
$E$ and $\tga^E$ is a section of $\Om^1(\ES)$\,.

More precisely, we first extend the operator
$\de^{-1}$ (\ref{del-1}) to $\Omb(\ES)$ and 
then define $\ga^E$ as a result of iterating 
the following equation (in degrees in $y$)  
\begin{equation}  
\label{iter-ga-E}
\ga^E = \G^E + \de^{-1} (\n \ga^E
+ A (\ga^E) + \frac{1}{2} [\ga^E, \ga^E])\,.
\end{equation}
Then $\ga^E$ satisfies the identity
\begin{equation}
\label{MC-ga}
D \ga^E + \frac{1}{2} [\ga^E , \ga^E] = 0\,,
\end{equation}
which immediately implies that $\tD^2=0$\,.

The differential $\tD$ (\ref{DDD-E1})
naturally extends to the sheaf of DGLAs
$\Omb(\Cbu(\ES))$ and to the sheaf of
DGLA modules $\Omb(\Cbd(\ES))$. Namely,
on $\Omb(\Cbu(\ES))$ $\tD$ is defined by
the formula
\begin{equation}
\label{DDD-coch}
\tD = D + [\pa \ga^E ,\, ]_G\,,
\end{equation}
and on  $\Omb(\Cbd(\ES))$ $\tD$ is
defined by
\begin{equation}
\label{DDD-ch}
\tD = D + R_{\pa\ga^E}\,,
\end{equation}
where $\ga^{E}$ is viewed locally as a section
of the sheaf $\Om^1(C^{-1}(\ES))$ and
$\pa$ denotes the Hochschild coboundary
operator.

Since the sheaf of DGLA modules
\begin{equation}
\begin{array}{c}
(\Omb (\Cbu(\SM)), D + \pa, [,]_{G}) \\[0.3cm]
 \downarrow_{\,mod} \\[0.3cm]
(\Omb (\Cbd(\SM)), D + \mb)
\end{array}
\label{para}
\end{equation}
is connected by a chain of quasi-isomorphisms 
to its cohomology (\ref{diag-T}), (\ref{diag-D}), 
(\ref{diag-K-Sh}) it suffices to show that 
(\ref{para}) is quasi-isomorphic to the 
sheaf of DGLA modules $\Cbd(\End(E))$ over 
the sheaf of DGLAs $\Cbu(\End(E))$. 
It is the sheaf of DGLA modules 
\begin{equation}
\begin{array}{c}
(\Omb (\Cbu(\ES)), \tD + \pa, [,]_{G}) \\[0.3cm]
 \downarrow_{\,mod} \\[0.3cm]
(\Omb (\Cbd(\ES)), \tD + \mb)
\end{array}
\label{para1}
\end{equation}
which allows us to do it. 

Indeed, generalizing the construction 
of the maps $\la_D$ and $\la_C$ in 
(\ref{diag-D}) we get the 
following embeddings
of the sheaves of DGLA modules
\begin{equation}
\begin{array}{ccc}
(\Omb(\Cbu(\ES)), \tD + \pa, [,]_{G}) &
\stackrel{\,\la^E_D}{\,\longleftarrow\,} & \Cbu(\End(E)) \\[0.3cm]
\downarrow^{R}_{\,mod}  & ~  &     \downarrow^{R}_{\,mod} \\[0.3cm]
(\Omb(\Cbd(\ES)), \tD + \mb) &\stackrel{\,\la^E_C}{\,\longleftarrow\,} 
&   \Cbd(\End(E))\,.
\end{array}
\label{diag-C-E}
\end{equation}
Similarly to propositions 7, 13, 15 in \cite{thesis} 
one can easily show that $\la^E_D$ and $\la^E_C$ are 
quasi-isomorphisms of the corresponding complexes 
of sheaves. 

Thus it remains to connect the sheaf of 
DGLA modules (\ref{para1}) to (\ref{para}) 
by a quasi-isomorphism. To do this we need the 
following auxiliary statement which is 
proved in a more general form in 
the concluding section 

\begin{lem}
\label{ona}
Let $a$, $b$, $c$, $d$ be
elements of a graded associative algebra with
the degrees
$$
\deg a = 0\,, \qquad \deg b = \deg c = \deg d = 1
$$
and let $a$ be nilpotent. 
If these elements satisfy the following
relations
\begin{equation}
\label{commut}
\begin{array}{ccc}
\displaystyle
[d , a] = b -\frac{c}{2}\,, & ~ & [b, a] = c\,, \\[0.3cm]
~ & [c,a]=0\,, & ~
\end{array}
\end{equation}
then
\begin{equation}
\label{resultat}
d \exp(a)  = \exp(a) (d  + b) \,. \qquad \Box
\end{equation}
\end{lem}

Let us pick a trivialization of $E$ over
a neighborhood $V$ of a point $p\in M$
and notice that on $V$ the initial Fedosov differential
(\ref{DDD}) on the sheaves $\Omb (\Cbu(\ES))$,
$\Omb (\Cbd(\ES))$ is well defined 
and the connection form $\ga^E$ (\ref{iter-ga-E})
can be viewed as a section of the sheaf
$\Om^1(\ES)\Big|_{V}$\,.

Furthermore, over $V$ the trace and
cotrace maps give the
following commutative diagram of
the quasi-isomorphisms of the 
sheaves of DGLAs and their modules
\begin{equation}
\begin{array}{ccc}
(\Omb (\Cbu(\SM)), D + \pa, [,]_{G})
&
\stackrel{\cotr}{\,\rightarrow\,}
&
(\Omb (\Cbu(\ES)), D + \pa, [,]_{G}) \\[0.3cm]
 \downarrow_{\,mod}
&
~
&
 \downarrow_{\,mod} \\[0.3cm]
(\Omb (\Cbd(\SM)), D+\mb) &
\stackrel{\tr}{\,\leftarrow\,}
&
(\Omb (\Cbd(\ES)), D + \mb)\,.
\end{array}
\label{tr-cotr}
\end{equation}

Now we notice that (\ref{MC-ga})
implies the identities
\begin{equation}
\label{D-b-R}
[(D + \mb), R_{\ga^E} ] =
R_{\pa \ga^E} - \frac{1}{2} R_{[\pa\ga^E,\ga^E]_G}\,,
\end{equation}
\begin{equation}
\label{D-b-br}
[(D + \pa), [\ga^E, \,]_G ] =
[\pa\ga^E, \,]_G -
 \frac{1}{2} [[\pa\ga^E,\ga^E]_G, \, ]_G\,,
\end{equation}
\begin{equation}
\label{Bianca}
[[\pa\ga^E,\ga^E]_G, \ga^E ]_G = 0\,,
\end{equation}
which allow us to apply lemma \ref{ona} 
to the algebras of operations on the 
sheaves $\Omb (\Cbu(\ES)$ and $\Omb (\Cbd(\ES)$\,.

Indeed, setting 
$$
a = [\ga^E, \, ]_G\,, \qquad  
b = [ \pa \ga^E , \, ]_G \,, \qquad
c = [ [\pa \ga^E , \ga^E]_G , \, ]_G \,, 
\qquad
d = D + \pa 
$$
and using (\ref{D-b-br}) and (\ref{Bianca})
we get that 
the map of complexes of sheaves (over $V$)
\begin{equation}
\label{twist-coch}
\exp (-[\ga^E\,, \, ]_G) : \Omb (\Cbu(\ES), D + \pa)
\to  (\Omb (\Cbu(\ES), \tD + \pa)
\end{equation}
is compatible with the corresponding 
differentials. 

Similarly, setting 
$$
a = R_{\ga^E}\,, \qquad  
b = R_{\pa \ga^E }  \,, \qquad
c = R_{[\pa \ga^E , \ga^E]_G }  \,, 
\qquad
d = D + \mb 
$$
and using (\ref{D-b-R}) and (\ref{Bianca})
we get that 
the map of complexes of sheaves (over $V$)
\begin{equation}
\label{twist-ch}
\exp (R_{\ga^E}) : (\Omb (\Cbd(\ES), \tD +\mb)
\to  (\Omb (\Cbd(\ES), D + \mb)
\end{equation}
is also compatible with the differentials.  

Combining these results with (\ref{tr-cotr})
we get the following commutative diagram 
of maps of sheaves of DGLAs and their modules
(over the neighborhood $V$)  
\begin{equation}
\begin{array}{ccc}
(\Omb (\Cbu(\SM)), D + \pa, [,]_{G})
&
\stackrel{\cotr^{tw}}{\,\rightarrow\,}
&
(\Omb (\Cbu(\ES)), \tD + \pa, [,]_{G}) \\[0.3cm]
 \downarrow_{\,mod}
&
~
&
 \downarrow_{\,mod} \\[0.3cm]
(\Omb (\Cbd(\SM)), D+\mb) &
\stackrel{\tr^{tw}}{\,\leftarrow\,}
&
(\Omb (\Cbd(\ES)), \tD + \mb)\,,
\end{array}
\label{tr-cotr-tw}
\end{equation}
where
\begin{equation}
\label{oni1}
\cotr^{tw} = \exp (-[\ga^E\,, \, ]_G) \circ \cotr\,,
\qquad
\tr^{tw} = \tr \circ \exp (R_{\ga^E})\,.
\end{equation}
Under changing the trivialization 
$\ga^E$ gets replaced by $\ga^E + \D$, 
where $\D$ is a one-form in $\Om^1(\ES)$ 
which does not involve the tangent coordinates 
$y^i$\,. Thus the maps (\ref{oni1})
do not depend on the choice of trivialization
because we deal normalized Hochschild 
(co)chains.  

Using the descending filtration associated 
to the exterior degree on the sheaves in 
(\ref{tr-cotr-tw}) one can easily show that 
the maps $\tr^{tw}$ and $\cotr^{tw}$  
are quasi-isomorphisms of complexes of sheaves.

The theorem is proved. $\Box$

\section{Concluding remarks}
In this section we show that lemma \ref{ona} admits 
an interesting generalization that involves the 
function 
\begin{equation}
\label{Todd}
f(x) = \frac{x}{e^x-1}
\end{equation}
from the definition of the Todd class. 

This generalization can be formulated as 
\begin{pred}
\label{ono}
Let $a$, $b$, $d$ be
elements of a graded associative algebra with
the degrees
$$
\deg a = 0\,, \qquad \deg b = \deg d = 1\,.
$$
Let $a$ be nilpotent and let 
$\al_1, \al_2, \dots, \al_k, \dots$ be
coefficients of the Taylor power series for 
the function
$$
\frac{x}{e^x -1} = 1 + \al_1 x + \al_2 x^2 + \dots 
$$
Then the relation 
\begin{equation}
\label{resultat1}
\exp(a) d  =  (d  + b) \exp(a) 
\end{equation}
holds if and only if 
\begin{equation}
\label{commut1}
\begin{array}{ccc}
\displaystyle
[a, d] = b + \sum_{k=1}^{\infty} \al_k 
(a d_a)^k\, b 
\end{array}
\end{equation}
\end{pred}
{\bf Proof.}  Since $a$ is a nilpotent element
all infinite power series in $a$ or 
in $a d_a$ are well defined.

Equation (\ref{commut1}) can be rewritten as
\begin{equation}
\label{vot}
[a, d] = f (a d_a)\, b\,,
\end{equation}
where $f=f(x)$ is given in (\ref{Todd})
and the operator $f(a d_a)$ is defined 
via the Taylor expansion of $f$ around 
the point $x=0$
$$
f(a d_a) = I d + \sum_{k=1}^{\infty} \al_k (a d_a)^k\,.
$$

It is clear that (\ref{vot}) holds 
if and only if
$$ 
b = g(a d_a) [a,d]\,,
$$
where $\displaystyle g(x) = \frac{e^x - 1}{x}$ 
and $g(a d_a)$ is also defined via the 
corresponding Taylor expansion. 

The latter equation is equivalent to 
$$
b = \exp (a d_a)\, d  - d    
$$
and this is exactly what we need to prove. $\Box$

Lemma \ref{ona} together with its generalization 
come to us as surprise. We hope that the relation 
with the formula for the Todd class 
can be helpful in deriving a version of the 
algebraic index theorem for the map  
(\ref{ind}).

~\\

\noindent\textsc{Department of Mathematics,
Northwestern University,
Evanston, IL 60208 \\
\emph{E-mail address:} {\bf vald@math.northwestern.edu}}

\end{document}